\documentclass[12pt, oneside]{amsart}   	% use "amsart" instead of "article" for AMSLaTeX format
\usepackage{geometry}
\usepackage{amsmath, amssymb}
\usepackage{epsfig}
\usepackage{epic,eepic}              		% See geometry.pdf to learn the layout options. There are lots.
\usepackage{graphicx}				% Use pdf, png, jpg, or eps? with pdflatex; use eps in DVI mode
\newtheorem{thm}{Theorem}[section]
\newtheorem{prop}[thm]{Proposition}
\newtheorem{lem}[thm]{Lemma}
\newtheorem{cor}[thm]{Corollary}

\theoremstyle{definition}
\newtheorem{definition}[thm]{Definition}

\newcommand\C{\mathbb{C}}
\newcommand\R{\mathbb{R}}
\newcommand\OO{\mathcal{O}}
\newcommand\M{\mathcal{M}}
\newcommand\F{\mathbb{F}}
\newcommand\PP{\mathbb{P}}

\numberwithin{equation}{section}

\title{Bounds of incidences between points and algebraic curves}

\begin{document}
\author{Hong Wang}
\address[H. Wang]{Department of Mathematics, Massachusetts Institute of Technology}
\email{hongwang@mit.edu}

\author{Ben Yang}
\address[B. Yang]{Department of Mathematics, Massachusetts Institute of Technology}
\email{beny@math.mit.edu}

\author{Ruixiang Zhang}
\address[R. Zhang]{Department of Mathematics, Princeton University}
\email{ruixiang@math.princeton.edu}

\begin{abstract}
We prove new bounds on the number of incidences between points and higher degree algebraic curves.  The key ingredient is an improved initial bound, which is valid for all fields. Then we apply the polynomial method to obtain global bounds on $\R$ and $\C$. 
\end{abstract}
\maketitle
\section{introduction}

The Szemer\'{e}di--Trotter theorem \cite{szemeredi1983extremal} says that for a finite set, $L$, of lines and a finite set, $P$, of points in $\mathbb{R}^2$, the number of incidences is less than a constant times $|P|^{\frac{2}{3}} |L|^{\frac{2}{3}} + |P| + |L|$. There have been several generalizations of this theorem. For example, Pach and Sharir \cite{pach1998number} allow simple curves that have $k$ degrees of freedom and multiplicity-type $C$: (i) for any $k$ distinct points there are at most $C$ curves in $L$ that pass through them, and (ii) any two distinct curves in $L$ have at most $C$ intersection points. This is summarized in the following theorem:

\begin{thm}[Pach--Sharir 98]\label{PSthm}
For a finite set $P$ of points in $\mathbb{R}^2$ and a finite set $L$ of simple curves which have $k$ degrees of freedom and multiplicity-type $C$. The number of incidences $|\mathcal{I}(P,L)|:=|\{(p,l)\in P\times L, p\in l\}|$ satisfies
\begin{equation}
|\mathcal{I}(P,L)|\lesssim_{C, k} |P|^{\tfrac{k}{2k-1}}|L|^{\tfrac{2k-2}{2k-1}}+|P|+|L|.
\end{equation}
\end{thm}
\subsection*{Notation} We use the asymptotic notation $X=\OO(Y)$ or $X\lesssim Y$ to denote the estimate $X\leq CY$ for some constant $C$. If we need the implicit constant $C$ to depend on additional parameters, then we indicate this by subscripts. For example, $X=\OO_{d}(Y)$ or $X\lesssim_{d} Y$ means that $X\leq C_{d}Y$ for some constant $C_{d}$ that  depends on $d$.

The main result of this paper is an improvement to Theorem \ref{PSthm} when $L$ is a set of higher degree algebraic curves. 

Let  $L$ be a finite set of algebraic curves of degree $\leq d$ in $\mathbb{R}^{2}$, any two of which do not share a common irreducible component. Let $P$ be a finite set of distinct points in $\mathbb{R}^{2}$. By B\'{e}zout's theorem, there is at most one curve in $L$ that goes through a subset of $P$ of size $d^{2}+1$. In the notation introduced by Pach and Sharir, a degree $d$ algebraic curve has $d^2+1$ degrees of freedom and multiplicity-type $d^2$. However, one may wonder whether $d^2 +1$ is a misleading definition of the degrees of freedom since generically $A:= {d+2 \choose 2} -1$ points determine a degree $d$ algebraic curve and $A\leq d^{2}+1$ when $d\geq 3$.

This suggests that Theorem ~\ref{PSthm} may still hold for degree $d$ curves with the  ``generic degree of freedom" $A$. Indeed, we prove that this is the case.

\begin{thm}\label{main theorem}
Let $d$ be a positive integer, $A={d+2\choose 2} -1$,  $L$ a finite set of degree $\leq d$ algebraic curves in $\mathbb{R}^{2}$ such that any two distinct algebraic curves do not share a common irreducible component, and $P$  a finite set of points in $\mathbb{R}^{2}$. Then,
\begin{equation}\label{estimatemainthm}
|\mathcal{I}(P,L)|\lesssim_{d} |P|^{\tfrac{A}{2A-1}}|L|^{\tfrac{2A-2}{2A-1}}+|P|+|L|.
\end{equation}
\end{thm}
This gives a better bound than Theorem ~\ref{PSthm} for degree $d$ algebraic curves when $d\geq 3$. It gives the same bound when $d=1$ or $d=2$.

We generalize Theorem ~\ref{main theorem} to algebraic curves parametrized by an algebraic variety. 
\begin{definition}
Given an integer $d\geq 1$ and a field $\F$, consider the Veronese embedding $\nu_{d}: \PP^{2}\rightarrow \PP^{A}$ given by:$$\nu_{d}: [x, y, z]\mapsto [x^{d},\dots, x^{i}y^{j}z^{k},\dots, z^{d}]_{i+j+k=d}.$$

We identify a degree $d$ curve in $\PP^{2}$ with the preimage of a hyperplane in $\PP^{A}$ or a point in the dual space $(\PP^{A})^{*}$. We say a subset $\M $ of the space of degree $d$ polynomials $\subset S_{d}$ is \emph{parametrized by an algebraic variety} $M$ if it is the preimage of $M\subset (\PP^{A})^{*}$, and we define the dimension of $\M$ to be $\dim M$.

\end{definition}

A consequence of the Theorem ~\ref{main theorem} is the following:

\begin{cor}\label{variety}
Given an integer $d\geq 1$,  $A={d+2\choose 2}-1$, and a subset $\mathcal{M}$ of the space of degree $d$ polynomials is parametrized by an algebraic variety $M$ of dimension $\leq k$. Let $P$ be a finite set of points in $\mathbb{RP}^2$, and $L$ a finite subset of $\M$ such that no two curves in $L$ share a common irreducible component. Then,
\begin{equation}
|\mathcal{I}(P,L)|\lesssim_{\mathcal{M}}|P|^{\tfrac{k}{2k-1}}|L|^{\tfrac{2k-2}{2k-1}}+|P|+|L|.
\end{equation}
\end{cor}

This generalization is helpful when we consider a family of curves with special properties, for example, parabolas, circles and a family of curves passing through a common point.

The proof of Theorem ~\ref{main theorem} is a standard application of the polynomial method (see for example \cite{dvir2009size} and \cite{guth2010erdos}) with the following new initial bound.  For an exposition of the polynomial method we refer the reader to \cite{kaplan2012simple}.

\begin{lem}[Initial bound]\label{trivial bound}
Let $\mathbb{F}$ be a field.  Given an integer $d\geq 1$, $A={d+2\choose 2} -1$, a finite set of algebraic curves, $L$,  in  $\mathbb{F}^{2}$ of degree $\leq d$ such that any two distinct curves in $L$ do not share a common component, and a finite set of points, $P$,  in $\mathbb{F}^{2}$. Then,
\begin{equation}
|\mathcal{I}(P,L)|\lesssim_{d} |P|^A+|L|.\label{initial bound}
\end{equation}
\end{lem}

Combining (\ref{initial bound}) with Solymosi and Tao's polynomial method and induction on the number of points \cite{solymosi2012incidence}, we obtain the following incidence theorem on complex space:

\begin{thm}\label{complex epsilon}
Given an integer $d\geq 1$,  $A={d+2\choose 2} -1$, a finite set of points $P$ in $\mathbb{C}^{2}$, a finite set of algebraic curves, $L$,  in $\mathbb{C}^{2}$ of degree $\leq d$ such that any two curves: (i) do not share a common irreducible component; (ii) intersect transversally at smooth points. Then, for any $\epsilon>0$,
\begin{equation}\label{estimateepsilonleft}
|\mathcal{I}(P,L)|\lesssim_{\epsilon,d} |P|^{\tfrac{A}{2A-1}+\epsilon}|L|^{\tfrac{2A-2}{2A-1}}+|P|+|L|
\end{equation}
and
\begin{equation}\label{estimateepsilonright}
|\mathcal{I}(P,L)|\lesssim_{\epsilon,d} |P|^{\tfrac{A}{2A-1}}|L|^{\tfrac{2A-2}{2A-1}+\epsilon}+|P|+|L|.
\end{equation}
\end{thm}
\subsection*{Acknowledgements} We would like to thank Larry Guth for encouraging us to work on the problem. He suggested that we try: improving the initial bound and working on algebraic subsets. We also thank Alex Townsend for reading drafts and for the suggestions.
\section{Proof of Theorem ~\ref{main theorem}}

%In this section we prove Theorem \ref{main theorem}. First we prove the ``initial bound'' Lemma \ref{trivial bound}.
%
\subsection{The initial bound}
In this subsection we prove the initial bound by double counting.
\begin{definition}
Given a point $p\in \PP^{2}$, let $H_{p}$ denote the corresponding hyperplane in $\PP^{A*}$ via the Veronese embedding and dual. Given a finite set of points $\Gamma=\{p_{1},\dots,p_{n}\}$, we define $m_{d}(\Gamma)=\dim(\cap H_{p_{i}})$, which characterizes the dimension of curves passing through all the points in $\Gamma$. In particular, if $m_{d}(\Gamma)=0$, there is at most one curve of degree $d$ that passes through $\Gamma$.
\end{definition}

\begin{proof}[Proof of Lemma \ref{trivial bound}]

We may remove curves containing fewer than $d^{2}+1$ points  by adding $\OO(|L|)$ to the bound. Now we assume that each curve contains more than $d^{2}+1$ points of $P$.

Fix a curve $l\in L$. We call an $A$-tuple $\Gamma'$ \emph{good}  if it is a subset of $\Gamma\subset l$ with $|\Gamma|=d^{2}+1$ and $m_{d}(\Gamma')=m_{d}(\Gamma)$. For any $(d^{2}+1)$-tuple $\Gamma\in l\cap P$, there exists a good $A$-tuple $\subset \Gamma$ since $L$ is parametrized by $\PP^{A}$.

 Since $\cap_{p\in \Gamma}H_{p}$ and $\cap_{p\in \Gamma'}H_{p}$ have the same dimension and $\cap_{p\in \Gamma}H_{p}\subseteq\cap_{p\in \Gamma'}H_{p}$, the two vector subspaces are in fact the same. In other words, any curve in $L$ passing through $\Gamma'$ must pass through all of $\Gamma$. Since curves in $L$ do not have a common irreducible component, by B\'{e}zout's theorem, every set of $d^{2}+1$ points determines a unique curve in $L$.  Hence, every good $A$-tuple determines a unique curve in $L$. There are at least ${|l\cap P| \choose d^2+1}/{ |l\cap P| - A\choose d^{2}+1-A}$ distinct good $A$-tuples $\Gamma'$ determining $l$. The number of good $A$-tuples determines the number of points in $l\cap P$ because ${|l\cap P| \choose d^2+1}/{ |l\cap P| - A\choose d^{2}+1-A}= \OO(|l\cap P|^{A})\gtrsim|l\cap P|$. On the other hand,  the number of $A$-tuples is ${|P|\choose A}= \OO(|P|^{A})$. Then,
$$|\mathcal{I}(P,L)| = \sum_{l \in L} |l\cap P|\lesssim |P|^{A}+|L|,$$
where the $|L|$ term comes from the first step when we deleted curves with fewer than $d^{2}+1$ points from $P$.
\end{proof}

\subsection{Polynomial method}\label{polynomialpartitioningsection}

Now we can apply the polynomial method to the initial bound and conclude the proof of Theorem \ref{main theorem}. We shall use the following polynomial partitioning proposition (see, for example, Theorem 4.1 in \cite{guth2010erdos}):

\begin{prop}\label{cell decomposition}
Let $P$ be a finite set of points in $\mathbb{R}^m$ and let $D$ be a positive integer. Then, there exists a nonzero polynomial $Q$ of degree at most $D$ and a decomposition
\[\mathbb{R}^{m}=\{Q=0\}\cup U_{1}\cup\cdots \cup U_{M}\]
into the hypersurface $\{Q=0\}$ and a collection $U_{1},\ldots ,U_{M}$ of open sets (which we call \emph{cells}) bounded by $\{Q=0\}$, such that $M =\OO_{m}( D^m)$ and that each cell $U_{i}$ contains $\OO_{m}(|P|/D^{m})$ points.
\end{prop}

\begin{proof}[Proof of Theorem \ref{main theorem}]
Applying Proposition ~\ref{cell decomposition}, we find a polynomial $Q$ of degree $D$ (to be chosen later) that partitions $\mathbb{R}^{2}$ into $M$ cells:
\[\mathbb{R}^{2}=\{Q=0\}\cup U_{1}\cup\cdots \cup U_{M},\]
where $M=\OO( D^2)$ and $U_{1},\ldots ,U_{M}$ are open sets bounded by  $\{Q=0\}$ and $P_{i}=U_{i}\cap P$. Let $L_{i}$ be the set of curves that have non-empty intersection with $U_{i}$. Then $|P_{i}|=\OO(|P|/D^{2})$. By B\'{e}zout's theorem we see every curve meets $\{Q=0\}$ at no more than $d\cdot D\leq \OO_{d}(D)$ components. Fix a curve $l\in L$, by Harnack's curve theorem~\cite{harnack1876ueber}, the curve itself  has $\OO_d (1)$ connected components. Moreover, since a component of $l$ is either contained in $U_{i}$ for some $i$ or it must meet the partition surface $\{Q=0\}$, $l$ can only meet $\OO_{d}(D)$ many $U_{i}$'s. Thus, we have the inequality $\sum |L_{i}|\leq D|L|$.

We may assume that every curve is irreducible, otherwise we can replace the curve with its irreducible components and the cardinality of curves increase by only a constant factor.
 Let $P_{cell}$ denote the points of $P$ in cells $U_{i}$ and let $P_{alg}$ denote those on the partition surface $\{Q=0\}$. Similarly, let $L_{alg}$ denote those curves that belong to $\{Q = 0 \}$ and let $L_{cell}$ be the union of the other curves. We deduce by Lemma \ref{trivial bound} and H\"{o}lder's inequality:
\begin{align}\label{polypartitionestimate}
|\mathcal{I}(P,L)|&= |\mathcal{I}(P_{cell}, L_{cell}) |+ |\mathcal{I}(P_{alg},L_{cell}) |+|\mathcal{I}(P_{alg}, L_{alg})|\nonumber\\
&\lesssim_{d}  \sum_i (|P_{i}|^{A}+|L_{i}|) + D|L_{cell}|+|\mathcal{I}(P_{alg}, L_{alg})|\nonumber\\
&\lesssim_{d}  |P|^{A}D^{-2(A-1)}+D|L|+|\mathcal{I}(P_{alg}, L_{alg})|.
\end{align}
In addition, we may assume that  $ |P|^{1/2}\leq |L|\leq |P|^{A}$, otherwise (\ref{estimatemainthm}) already holds either by Lemma \ref{trivial bound} or by another initial bound $|\mathcal{I} (P, L)| \lesssim |P| + |L|^2$ (every two curves intersect at at most $\OO(1)$ points). In this case, we may choose $D =\OO_{d}( |P|^{\tfrac{A}{2A-1}} |L|^{-\tfrac{1}{2A-1}})$ and $D\leq|L|/2$. Then, the first two terms on the right-hand side of (\ref{polypartitionestimate}) are $O_{d}(  |P|^{\tfrac{A}{2A-1}}|L|^{\tfrac{2A-2}{2A-1}})$. Since $|L_{alg}| \leq D \leq \frac{|L|}{2}$, we can perform a dyadic induction on $|L|$. By repeating the above process  on $L_{alg}$, we obtain the following:
\begin{align}
|\mathcal{I}(P,L)|&\lesssim_{d} \sum_{ 0\leq i \leq \log(|L|/|P|^{1/2})} |P|^{\tfrac{A}{2A-1}}(2^{-i}|L|)^{\tfrac{2A-2}{2A-1}} + |P|+|L|\nonumber\\&\lesssim_{d}|P|^{\tfrac{A}{2A-1}}|L|^{\tfrac{2A-2}{2A-1}}+|P|+|L|.\nonumber
\end{align}
 The induction stops when $i> \log(|L|/|P|^{1/2})$ because when $2^{-i}L\lesssim |P|^{1/2}$ the number of incidences in the $i$-th step is bounded by $\OO( |P| )$. This proves (\ref{estimatemainthm}).

\end{proof}

\section{An estimate for parametrized curves}

We prove an initial bound with parametrized curves, which implies Corollary \ref{variety}.
We first state two propositions that we will need to prove the initial bound.
\begin{prop}(see \cite{hartshorne1977algebraic}, Ch I, Exercise 1.8)\label{dimension diminute}
If $V$ is an $r$-dimensional variety in $\mathbb{F}^{d}$, and $P :\mathbb{F}^{d}\rightarrow V$ is a polynomial which is not identically zero on $V$, then every component of $V\cap \{P=0\}$ has dimension $r-1$.
\end{prop}

\begin{prop}(see \cite{fulton1984intersection}, Section 2.3)\label{refine Bezout}
Let $V_{1},\ldots V_{s}$ be subvarieties of $\mathbb{FP}^{N}$, and let $Z_{1},\ldots ,Z_{r}$ be the irreducible components of $V_{1}\cap\cdots \cap V_{r}$. Then,
\[\sum_{i=1}^{r} \deg(Z_{i})\leq \prod_{j=1}^{s} \deg(V_{j}).\]
\end{prop}
\begin{lem}\label{algebraicsettrivialbound}
With the same setting and notation as in Corollary \ref{variety}, we have
\begin{equation}
|\mathcal{I}(P,L)|\lesssim_{\mathcal{M}}|P|^{k}+|L|.
\end{equation}
\end{lem}

\begin{proof}
Without loss of generality we assume that $\mathbb{F}$ is an algebraically closed field, every curve in $L$ is irreducible, of degree $d$, and contains more than $d^{2}+1$ points from $P$. A curve $l\in L$ corresponds to a point of intersection  $\cap_{p\in l\cap P}H_{p}\cap M$. Comparing with the proof of Lemma \ref{trivial bound}, it suffices to prove that for every $(d^{2}+1)$-tuple $\Gamma$ in $l\cap P$, there is a $\Gamma' \subseteq \Gamma$ such that  $\cap_{p\in \Gamma'} H_{p}\cap M$ contains at most $O_{\mathcal{M}} (1)$ curves and $|\Gamma'|=k$. This follows from Proposition~\ref{dimension diminute} and Proposition~\ref{refine Bezout} above.

Indeed, by iterately applying Proposition ~\ref{dimension diminute}, we can choose $\Gamma'$ such that $| \Gamma'|=k$ and $\cap_{p\in\Gamma'}H_{p}\cap\mathcal{M}$ has dimension $0$. By Proposition \ref{refine Bezout}, the cardinality of $\cap_{p\in\Gamma'}H_{p}\cap\mathcal{M}$ is bounded by a constant depending on $k$ and $\deg \mathcal{M}$.

\end{proof}

\section{A theorem on the complex field with $\epsilon$}
In this section, we follow the approach of \cite{solymosi2012incidence} to sketch a proof of Theorem \ref{complex epsilon}. The idea is to partition the point set $P$ with a polynomial of constant degree and then use induction on the size of $P$. In other words, the degree does not depend on the size of $P$ and $L$. With an $\epsilon$ loss in the exponents, one can perform induction on $|P|$, which controls incidences in the complement of the partition surface (in the cells). Here, a constant degree implies constant complexity, and we can use dimension reduction to estimate incidences on the partition surface.  

\begin{proof}[Proof of Theorem \ref{complex epsilon}]

In our proof, $C$ is a constant depending only on $d$ and $\epsilon$ which may vary from place to place. $C_0, C_1$ and $C_2$ are positive constants to be chosen later, where $C_{0}, C_1 >2$ are sufficiently large depending on $d$ and $\epsilon$, and $C_{2}$ is sufficiently large depending on $C_{1}$, $C_{0}$, $d$ and $\epsilon$.

 We do an induction on $|P|$. Suppose for $|P'|\leq \tfrac{|P|}{2}$ and $|L'| \leq |L|$, we have by the induction hypothesis,
\begin{equation}\label{hypothesis}
|\mathcal{I}(P',L')|\leq C_{2}|P'|^{\tfrac{A}{2A-1}+\epsilon}|L'|^{\tfrac{2A-2}{2A-1}}+C_{0}(|P'|+|L'|).
\end{equation}
Our goal is to prove
\begin{equation}\label{inductiveclaim}
|\mathcal{I}(P,L)|\leq C_{2}|P|^{\tfrac{A}{2A-1}+\epsilon}|L|^{\tfrac{2A-2}{2A-1}}+C_{0}(|P|+|L|).
\end{equation}
We apply Proposition \ref{cell decomposition} to $D = C_{1}$ on $\mathbb{C}^2 \simeq \mathbb{R}^{4}$ and obtain the following partition:
\begin{equation}
\mathbb{R}^{4}=\{Q=0\}\cup U_{1}\cup\cdots \cup U_{M}.
\end{equation}

Here, $Q: \mathbb{R}^{4}\rightarrow\mathbb{R}$ has degree at most $C_{1}$, $M =\OO( C_{1}^{4})$, and $|P_{i}| = |P \cap U_{i}| = O(\tfrac{|P|}{C_{1}^{4}})\leq \tfrac{|P|}{2}$. We denote $L_{i}$ to be the set of curves in $L$ with nonempty intersection with $U_{i}$. Thus, by the induction hypothesis we have,
\begin{align}
|\mathcal{I}(P_{i},L_{i})| \leq & C_{2}|P_{i}|^{\tfrac{A}{2A-1}+\epsilon}| L_{i}|^{\tfrac{2A-2}{2A-1}}+C_{0}(|P_{i}| + |L_{i}|)\nonumber\\
\leq & C [C_{2}C_{1}^{-4(\tfrac{A}{2A-1}+\epsilon)}|P|^{\tfrac{A}{2A-1}+\epsilon}|L_{i}|^{\tfrac{2A-2}{2A-1}}+C_{0}(\tfrac{|P|}{C_1^4}+|L_{i}|)].
\end{align}

For $l$ belonging to some $L_{i}$, we apply a result in real algebraic geometry that implies the number of connected components of $l\setminus \{Q=0\}$ is at most $\OO_{d} (C_{1}^{2})$ (see Theorem A.2 of \cite{solymosi2012incidence}, \cite{milnor1964betti},\cite{petrovskii1949topology},\cite{thom1965homologie}). We deduce that
\begin{equation}
\sum_{i=1}^{M}|L_{i}| \leq CC_{1}^{2}|L|.
\end{equation}

Adding up $|\mathcal{I}(P_i,L_{i})|$ and applying H\"{o}lder's inequality, we obtain
\begin{align}
|\mathcal{I}(P_{cell},L_{cell})| = & \sum_{i=1}^{M}|\mathcal{I}(P_{i},L_{i})|\nonumber\\
\leq & C(C_{2}C_{1}^{-4(\tfrac{A}{2A-1}+\epsilon)}|P|^{\tfrac{A}{2A-1}+\epsilon}(\sum_{i=1}^{M}|L_{i}|^{\tfrac{2A-2}{2A-1}})+C_0 (|P|+C_{1}^{2}|L|))\nonumber\\
\leq & C(C_{1}^{-4\epsilon}C_{2}|P|^{\tfrac{A}{2A-1}+\epsilon}|L|^{\tfrac{2A-2}{2A-1}}+C_{0}(|P|+C_{1}^{2}|L|)).
\end{align}

Now we recall the two trivial bounds (\ref{L2+P}) explained in the case of real space and by Lemma~\ref{trivial bound}, we have:

\begin{equation}\label{L2+P}
|\mathcal{I}(P,L)|\lesssim_{d} |L|^{2}+|P|,~~~~~
|\mathcal{I}(P,L)|\lesssim_{d} |P|^{A}+|L|.
\end{equation}
Thus, one may assume that $|P|^{\tfrac{1}{2}}\lesssim_{d}|L|\lesssim_{d}|P|^{A}$, otherwise we immediately have $|\mathcal{I}(P,L)|\lesssim_{d}|P|+|L|$ and it suffices to choose $C_0$ larger than the implicit constant. With this assumption, we have
\begin{equation}\label{PcellLcell}
|\mathcal{I}(P_{cell},L_{cell})| \leq  C(C_{1}^{-4\epsilon}C_{2}+C_{0}(|P|^{-\epsilon}+C_{1}^{2}|P|^{-\epsilon}))|P|^{\tfrac{A}{2A-2}+\epsilon}|L|^{\tfrac{2A-2}{2A-1}}.
\end{equation}
If the following inequality is given
\begin{equation}\label{PalgL}
\mathcal{I}(P_{alg}, L)\lesssim_{C_1}|P|^{\tfrac{A}{2A-1}}|L|^{\tfrac{2A-2}{2A-1}} + |P|+|L|,
\end{equation}
then $ \lesssim_{C_1} |P|^{\tfrac{A}{2A-1}}|L|^{\tfrac{2A-2}{2A-1}}$ by assumption. Hence, we can combine this with (\ref{PcellLcell}) and a careful choice of $C_0, C_1$ and $C_2$ gives us (\ref{inductiveclaim}). When $|P|  = 1$, (\ref{inductiveclaim}) is trivial, and we obtain (\ref{estimateepsilonleft}). For (\ref{estimateepsilonright}) the argument is similar and is omitted.

Finally, (\ref{PalgL}) follows from Proposition ~\ref{inductionondim} below when $r=3$, $D=C_{1}$ and $\Sigma=\{Q=0\}$.
\end{proof}

\begin{prop}\label{inductionondim}
Let $P$ and $L$ be as in Theorem ~\ref{complex epsilon}, $0\leq r\leq 3$, and let $\Sigma$ be a subvariety in $\mathbb{C}^2 \simeq \mathbb{R}^{4}$ of (real) dimension $\leq r$ and of degree $\leq D$. Then,
\begin{equation}
\mathcal{I}(P\cap \Sigma, L) \lesssim_{D} |P|^{\tfrac{A}{2A-1}}|L|^{\tfrac{2A-2}{2A-1}}+|P|+|L|.
\end{equation}
\end{prop}

\begin{proof}
When $r=0$, $\Sigma$ is a single point and the inequality trivially holds.

When $r=1$, we decompose $\Sigma=\Sigma_{1}\cup\Sigma_{2}$, where every component  of $\Sigma_{1}$ belongs to some curve in $L$ and $\Sigma_{2}$ do not have a common component with curves in $L$. Then, $|\mathcal{I}(P\cap \Sigma_{1})|\lesssim_{D}|P|$ and $|\mathcal{I}(P\cap \Sigma_{2})|\lesssim_{D, d}|L|$. 

Now we deal with the case $r=2$. By an algebraic geometry result (see, for example, Corollary 4.5 in \cite{solymosi2012incidence}), one can decompose $\Sigma$ into smooth points on subvarieties
\[\Sigma = \Sigma^{smooth}\cup\Sigma_{i}^{smooth},\]
where $\Sigma_{i}$'s are subvarieties of $\Sigma$ of dimension $\leq 1$ and of degree $\OO_{D}(1)$. The number of $\Sigma_{i}$'s is at most $\OO_{D}(1)$. It suffices to bound $\mathcal{I}(P\cap\Sigma^{smooth}, L)$.

If $l_{1}, l_{2} \in \Sigma$ intersect at $p\in \Sigma^{smooth}$, by considering the tangent space and the transverse assumption, we know that $p$ is a singular point of $l_{1}$ or $l_{2}$. Since each curve has $\OO_d (1)$ singular points, we obtain
\begin{equation}\label{SigmaLalg}
\mathcal{I}(P\cap\Sigma^{smooth}, L_{alg}) \leq |P| + \OO_d (1)|L|.
\end{equation}

It remains to estimate $\mathcal{I}(P\cap\Sigma^{smooth}, L_{cell})$. If $l$ does not belong to $\Sigma$, then by Corollary 4.5 of \cite{solymosi2012incidence} the intersection of $l$ and $\Sigma$ can be decomposed as $l\cap\Sigma=\cup_{j=0}^{J (l)}l_{j}$ for some $J(l) \leq \OO_{D}(1)$, where $l_{j}$ is an algebraic variety of dimension $\leq 1$ and of degree $\OO_D (1)$  for each $1\leq j\leq J(l)$. Let $\mathcal{I}_{l, j}$ denote the set $\{p \in P: p \in l_{j}\}$, we obtain
\begin{equation}
|\mathcal{I}(P\cap \Sigma^{smooth}, L_{cell})|\leq \sum_{l, j: j \leq J(l)} |\mathcal{I}_{l, j}|.
\end{equation}

If $l_{j}$ is not the union of $\OO_D (1)$ points, then $l_{j}$ belongs to a unique $l$ because distinct curves in $L$ do not share a common component. By taking a generic projection from $\mathbb{R}^4$ to $\mathbb{R}^2$, we can apply arguments in the proof of Theorem ~\ref{main theorem}: use the initial bound given by $L$ and $P$, then apply the polynomial method.

When $r=3$,  we can repeat  the proof of $r=2$ assuming that the bound holds for $r\leq 2$.

\end{proof}

We also have the following corollary for complex curves parametrized by an algebraic variety:
\begin{cor}\label{complex variety}
Given a finite point set $P\in \mathbb{C}^{2}$, an integer $d\geq 1$, $A={d+2\choose 2} -1$ and a subset $\mathcal{M}\in S_d$ parametrized by an algebraic variety of dimension $\leq k$. Let $L$ be a finite subset of $\mathcal{M}$ such that any two distinct curves of $L$ do not share a common component and intersect transversally at smooth points. Then, for any sufficiently small $\epsilon>0$,
\begin{equation}\label{estimatecomplexepsilonleft}
|\mathcal{I}(P,L)|\lesssim_{\epsilon, \mathcal{M}} |P|^{\tfrac{k}{2k-1}+\epsilon}|L|^{\tfrac{2k-2}{2k-1}}+|P|+|L|
\end{equation}
and
\begin{equation}\label{estimatecomplexepsilonright}
|\mathcal{I}(P,L)|\lesssim_{\epsilon, \mathcal{M}} |P|^{\tfrac{k}{2k-1}}|L|^{\tfrac{2k-2}{2k-1}+\epsilon}+|P|+|L|.
\end{equation}
\end{cor}
\bibliographystyle{amsalpha}
\bibliography{stageref}

\providecommand{\bysame}{\leavevmode\hbox to3em{\hrulefill}\thinspace}
\providecommand{\MR}{\relax\ifhmode\unskip\space\fi MR }
\providecommand{\MRhref}[2]{%
  \href{http://www.ams.org/mathscinet-getitem?mr=#1}{#2}
}
\providecommand{\href}[2]{#2}
\begin{thebibliography}{KMS12}

\bibitem[Dvi09]{dvir2009size}
Zeev Dvir, \emph{On the size of kakeya sets in finite fields}, Journal of the
  American Mathematical Society \textbf{22} (2009), no.~4, 1093--1097.

\bibitem[Ful84]{fulton1984intersection}
William Fulton, \emph{Intersection theory}, vol. 1998, Springer Berlin, 1984.

\bibitem[GK10]{guth2010erdos}
Larry Guth and Nets~Hawk Katz, \emph{On the erdos distinct distance problem in
  the plane}, arXiv preprint arXiv:1011.4105 (2010).

\bibitem[Har76]{harnack1876ueber}
Axel Harnack, \emph{Ueber die vieltheiligkeit der ebenen algebraischen curven},
  Mathematische Annalen \textbf{10} (1876), no.~2, 189--198.

\bibitem[Har77]{hartshorne1977algebraic}
Robin Hartshorne, \emph{Algebraic geometry}, vol.~52, springer Verlag, 1977.

\bibitem[KMS12]{kaplan2012simple}
Haim Kaplan, Ji{\v{r}}{\'\i} Matou{\v{s}}ek, and Micha Sharir, \emph{Simple
  proofs of classical theorems in discrete geometry via the guth--katz
  polynomial partitioning technique}, Discrete \& Computational Geometry
  \textbf{48} (2012), no.~3, 499--517.

\bibitem[Mil64]{milnor1964betti}
John Milnor, \emph{On the betti numbers of real varieties}, Proceedings of the
  American Mathematical Society \textbf{15} (1964), no.~2, 275--280.

\bibitem[PO49]{petrovskii1949topology}
Ivan~Georgievich Petrovskii and Olga~Arsen'evna Oleinik, \emph{On the topology
  of real algebraic surfaces}, Izvestiya Rossiiskoi Akademii Nauk. Seriya
  Matematicheskaya \textbf{13} (1949), no.~5, 389--402.

\bibitem[PS98]{pach1998number}
J{\'a}nos Pach and Micha Sharir, \emph{On the number of incidences between
  points and curves}, Combinatorics, Probability and Computing \textbf{7}
  (1998), no.~01, 121--127.

\bibitem[ST12]{solymosi2012incidence}
J{\'o}zsef Solymosi and Terence Tao, \emph{An incidence theorem in high er
  dimensions}, Discrete \& Computational Geometry \textbf{48} (2012), no.~2,
  255--280.

\bibitem[STJ83]{szemeredi1983extremal}
Endre Szemer{\'e}di and William~T Trotter~Jr, \emph{Extremal problems in
  discrete geometry}, Combinatorica \textbf{3} (1983), no.~3-4, 381--392.

\bibitem[Tho65]{thom1965homologie}
Ren{\'e} Thom, \emph{Sur l¡¯homologie des vari{\'e}t{\'e}s alg{\'e}briques
  r{\'e}elles}, Differential and combinatorial topology (1965), 255--265.

\end{thebibliography}

\end{document}